# TREE REPRESENTATIONS OF GALOIS GROUPS

Nigel Boston

University of Illinois at Urbana-Champaign

**Introduction.**

Much work has gone into matrix representations of Galois groups, but there is a whole new class of naturally occurring representations that have as yet gone almost unnoticed. In fact, it is well-known in various areas of mathematics that the main sources of totally disconnected groups are matrix groups over local fields AND automorphism groups of locally finite trees [12]. It is perhaps surprising then that representations of Galois groups into the latter have been almost ignored, while at the same time Galois representations into the former have been enormously effective in resolving long-standing problems in number theory.

These "tree" representations are important as regards topics such as the unramified Fontaine-Mazur conjecture [6]. This conjecture states that any $p$-adic representation of the Galois group $G$ of an extension unramified at $p$ (and ramified at only finitely many primes) should have finite image. In other words, $p$-adic representations say little about such Galois groups. This paper proposes the conjecture that these Galois groups should, on the other hand, have representations with large image (measured by Hausdorff dimension) in the automorphism group of a rooted tree. Thus, they might allow us to investigate the structure of the Galois group of infinite pro-$p$ extensions (such as Hilbert $p$-class towers), something unapproachable by standard $p$-adic representation methods.

**The Conjecture.**

If $G$ is any infinite, finitely generated pro-$p$ group, then it has minimally infinite quotients called *just-infinite* pro-$p$ groups [8]. If $K$ is any number field and $S$ is a finite set of places of $K$, then we denote by $G_{K,S}$ the Galois group over $K$ of a maximal extension unramified outside $S$. If $p$ is a rational prime and $S$ contains no places above $p$, then the conjecture of Fontaine and Mazur says that any finite-dimensional $p$-adic representation of $G_{K,S}$ should have finite image. In a previous paper [2], I explained how this is equivalent to the conjecture that the just-infinite pro-$p$ quotients of all such $G_{K,S}$ are not linear over $\mathbf{Z}_p$. The question then arises as to just what structure these just-infinite groups do have. For convenience, we call these groups *uji groups* below (standing for unramified just-infinite).

The author was partially supported by NSF grant DMS 99-70184.

The attempted classification [10] of just-infinite pro-$p$ groups has so far yielded four classes of just-infinite pro-$p$ groups, namely

I. Solvable (and so linear over $\mathbf{Z}_p$).
II. Nonsolvable and linear over $\mathbf{Z}_p$.
III. Nonsolvable and linear over $\mathbf{F}_p[[T]]$.
IV. The rest!

At present, the rest consists of certain subgroups of the Nottingham groups $R_k$ ($k$ a finite field) (in particular open subgroups and the Fesenko groups [5]), where $R_k$ consists of the automorphisms $T \mapsto T + a_2 T^2 + a_3 T^3 + ...$ of $k[[T]]$, and Grigorchuk-type (or 'branch') groups [7], which are certain pro-$p$ automorphism groups of the $p$-ary rooted tree. This tree has $p^n$ vertices at distance $n$ from its root with each vertex above the root having valency $p + 1$. In fact, Grigorchuk [7] has recently proved that *every* just-infinite pro-$p$ group either is branch or contains an open subgroup of the form $H \times ... \times H$ (finitely many factors), where $H$ is hereditarily just-infinite (i.e. every open subgroup of $H$ is just-infinite).

Uji groups cannot be of type I by class field theory and the conjecture of Fontaine-Mazur implies they are not of type II. In [2] I gave evidence for a conjecture that implies that the uji groups cannot be of type III either. In this note, we explore the possibility that all uji groups lie on the same side of Grigorchuk's dichotomy, namely are branch. This is a strengthening and an elaboration of the unramified Fontaine-Mazur conjecture.

Note also that each of the types of just-infinite pro-$p$ groups is attached to a kind of Galois representation. Linear representations over $\mathbf{Z}_p$ and $\mathbf{F}_p[[T]]$ are well-known and important. The field of norms construction yields representations of local Galois groups into $R_k$, something that Fesenko has attempted to globalize. How about Galois representations into branch groups?

At first glance there are too many of these to be useful since every finitely generated pro-$p$ group embeds in the pro-$p$ automorphism group $W$ of the $p$-ary rooted tree $T$. Our conjecture, however, concerns representations with large image. This is a familiar idea from the theory of $p$-adic representations, although here 'large' cannot mean 'of finite index' since $W$ is not finitely generated, whereas our pro-$p$ Galois groups are. To define 'large', we need the notion of Hausdorff dimension [1].

Let $W_n$ be the quotient of $W$ given by its action on the subtree of vertices of distance $\leq n$ from the root of $T$. If $G$ is a closed subgroup of $W$, its Hausdorff dimension is defined to be $\liminf_{n \to \infty} \frac{\log |G_n|}{\log |W_n|}$, where $G_n$ is the image of $G$ in $W_n$. So, for instance, the first branch group discovered (by Grigorchuk) [7] has Hausdorff dimension $5/8$. If $G$ is any pro-$p$ group, then we define its Hausdorff dimension to be the supremum of its Hausdorff dimensions over all embeddings into $W$.

**Conjecture 1.** A just-infinite pro-$p$ group is branch if and only if its Hausdorff dimension is nonzero.

**Conjecture 2.** The just-infinite pro-$p$ quotients of $G_{K,S}$ are branch.

**Corollary to Conjecture 2.** The unramified Fontaine-Mazur conjecture [6] (and its generalization in [2]) holds.

Conjecture 1 is purely group-theoretical. Combining the two conjectures we see that when $G_{K,S}$ has infinite pro-$p$ quotients, it should have maps to $W$ with 'large' image in the sense of nonzero Hausdorff dimension. These arise by mapping $G_{K,S}$ onto a just-infinite quotient $J$ and then embedding $J$ in $W$ with nonzero Hausdorff dimension.

**Further Discussion.**

We concentrate on the $p = 2$ case here. There are straightforward generalizations to other $p$, but this is computationally most convenient. Let $W_1 = C_2$ and $W_n = W_{n-1} \wr C_2$ ($n \geq 2$). Then $W_n$ is a finite 2-group of order $2^{2^n - 1}$, that naturally arises as the automorphism group of the binary tree $T_n$, which consists of a single root vertex with two edges above each vertex to the next level extending to level $n$ (in fact, $W_n$ is isomorphic to the Sylow 2-subgroup of the symmetric group on $2^n$ letters). Thus $W_{n-1}$ is naturally a quotient of $W_n$ and we set $W = \varprojlim W_n$, an infinitely generated pro-2 group, which is the pro-2 automorphism group of the infinite binary tree $T$. We will be interested in representations of $G_{\mathbf{Q}} := \text{Gal}(\overline{\mathbf{Q}}/\mathbf{Q})$ into $W = \text{Aut}(T)$.

The only place where representations of this kind have been studied so far is by Odoni [9], who in work completed by Stoll [11], showed that the Galois group of the $n$th iterate of $x^2 + 1$ is $W_n$, so that letting $n \to \infty$, we obtain a surjective representation $G_{\mathbf{Q}} \to W$. In showing surjectivity, their aim was to prove that new primes were ramified at each level.

An example of a representation $G_{\mathbf{Q},S} \to W_4$ with large image and $S$ not containing 2 is given in the last section of my work with Leedham-Green [3].

Some questions that arise from this are:

(i) Is there a natural way to construct tree representations, perhaps by using a series of quadratic covering maps (the Odoni-Stoll case corresponding to the case of $\mathbf{P}^1 \to \mathbf{P}^1$ given by $x \mapsto x^2 + 1$)? Algebraic geometry could then reenter the study of $G_{K,S}$.

(ii) Given a tree representation ramified at finitely many places, do the images of Frobenius elements carry arithmetical information? Are there objects related to the representation in a way similar to modular forms being related to $p$-adic representations via traces of Frobenius? Note that the conjugacy classes of the automorphism groups of $p$-ary trees have received much study lately - see e.g. [4].

## BIBLIOGRAPHY

1. Y.Barnea and A.Shalev, Hausdorff dimension, pro-$p$ groups, and Kac-Moody algebras, Trans. Amer. Math. Soc. **349** (1997), no. 12, 5073–5091.

2. N.Boston, Some Cases of the Fontaine-Mazur Conjecture II, J. Number Theory **75** (1999), 161–169.

3. N.Boston and C.R.Leedham-Green, Explicit computation of Galois $p$-groups unramified at $p$ (preprint).

4. A.M.Brunner and S.Sidki, On the automorphism group of the one-rooted binary tree, J. Algebra **195** (1997), no. 2, 465–486.

5. I.Fesenko, On just infinite pro-$p$-groups and arithmetically profinite extensions of local fields, J. Reine Angew. Math. **517** (1999), 61–80.

6. J.-M.Fontaine and B.Mazur, Geometric Galois representations, *in* "Elliptic curves and modular forms, Proceedings of a conference held in Hong Kong, December 18-21, 1993," International Press, Cambridge, MA and Hong Kong.

7. R.Grigorchuk, Article in "New Horizons in pro-$p$ Groups" (eds. du Sautoy, Segal, Shalev), Birkhauser, Boston 2000.

8. G.Klaas, C.R.Leedham-Green, and W.Plesken, Linear pro-$p$ groups of finite width, Springer Lecture Notes in Math 1674, 1997.

9. R.W.K.Odoni, Realising wreath products of cyclic groups as Galois groups, Mathematika **35** (1988), no. 1, 101–113.

10. A.Shalev, On almost fixed point free automorphisms, J. Algebra **157**, (1993), 271-282.

11. M.Stoll, Galois groups over **Q** of some iterated polynomials, Arch. Math. (Basel) **59** (1992), no. 3, 239–244.

12. G.Willis, Totally disconnected, nilpotent, locally compact groups, Bull. Austral. Math. Soc. **55** (1997), no. 1, 143–146a.

Department of Mathematics, University of Illinois, Urbana, Illinois 61801
*E-mail address*: `boston@math.uiuc.edu`